\documentclass{article}
\usepackage{latexsym,amsmath, amssymb,graphicx,tikz}
\newcommand{\new}{\newcommand}

\new{\bracket}[1]{\left\langle #1 \right \rangle}
\new{\parens}[1]{\!\left( #1 \right)}
\new{\sbrace}[1]{\!\left[ #1 \right]}
\new{\cbrace}[1]{\!\left\{ #1 \right\}}
\new{\abs}[1]{\left| #1 \r|}
\numberwithin{equation}{section}

\begin{document}
\nocite{ST67}

\title{
 South Pointing Chariot: An Invitation to Differential Geometry
}
\author{Stephen Sawin\\Fairfield University}
\maketitle

\begin{abstract}
We introduce the south-pointing chariot, an intriguing mechanical
device from ancient China.  We use its ability to keep track of a
global direction as it travels on an arbitrary path as a tool to explore
the geometry of curved surfaces.  This takes us  as far as a famous result of
Gauss on the impossibility of a faithful map of the globe, which
started off the field of differential geometry.  The reader should get
a view into how geometers think and an introduction to important early
results in the field, but should need no more than a solid background in calculus
(ideally through multivariable calculus).  This is achieved  by relying
on the reader's visual intuition.  
\end{abstract}

\section{Some History, and Some Engineering}

Legend tells of  Huangdi, the Yellow
Emperor, tamer of beasts real and supernatural, father of the Huaxia
people who formed what would become China.   Also  inventor of the bow sling,
the Chinese calendar, early astronomy, a zither, a Chinese version of
football, and depending on whom you ask early Chinese writing and the growing of
cereal.  Four thousand years ago with the aid of his many tamed beasts he fought a famous battle with his arch-nemesis, the
bronze-headed Chi-You and his $81$ horned and four-eyed brothers. Chi You confounded the Yellow
Emperor by breathing out a  thick fog, but the Emperor invented,
apparently on the spur of the moment, an ingenious ``south pointing
chariot'' with a figure on the top connected to the wheels in such a
fashion that it always pointed south.  This contraption led his army out of the
fog and to victory.

It may not suprise you that modern scholars do not give complete
credence to the details
of this story. Among other things it is doubtful that the technology
to build such a chariot could have been achieved so early.   This skepticism was shared even by some early Chinese
figures, such as Permanent Counsellor Gaotang Long and 
General Qin Lang during the Three Kingdoms period in the third century
CE.  Their contemporary, the  engineer Ma Jun, replied to their doubt cuttingly with ``Empty
arguments with words cannot (in any way) compare with a test which
will show practical results'' and silenced them by inventing and
building a version of Huangdi's chariot.  No
designs or physical evidence remains, but scholars do generally
credit that he built such a thing.  That's no tamed dragon, but still
quite remarkable.  It seems to require gearing unknown until 1720 CE.  The south
pointing chariot was lost and reinvented several times in Chinese
history, and while there are reports of its use in military and
navigation applications, all we know it was used for is ceremonial
processions \cite{Needham65}.

\section{How it Works}
 
If you search on a term such as ``south pointing chariot'' you will
find videos of recreations that will illustrate how the gears allow
this magic to happen.  We'll try here with words and a few pictures.
First the differential. 
\begin{figure}
\includegraphics[width=5cm]{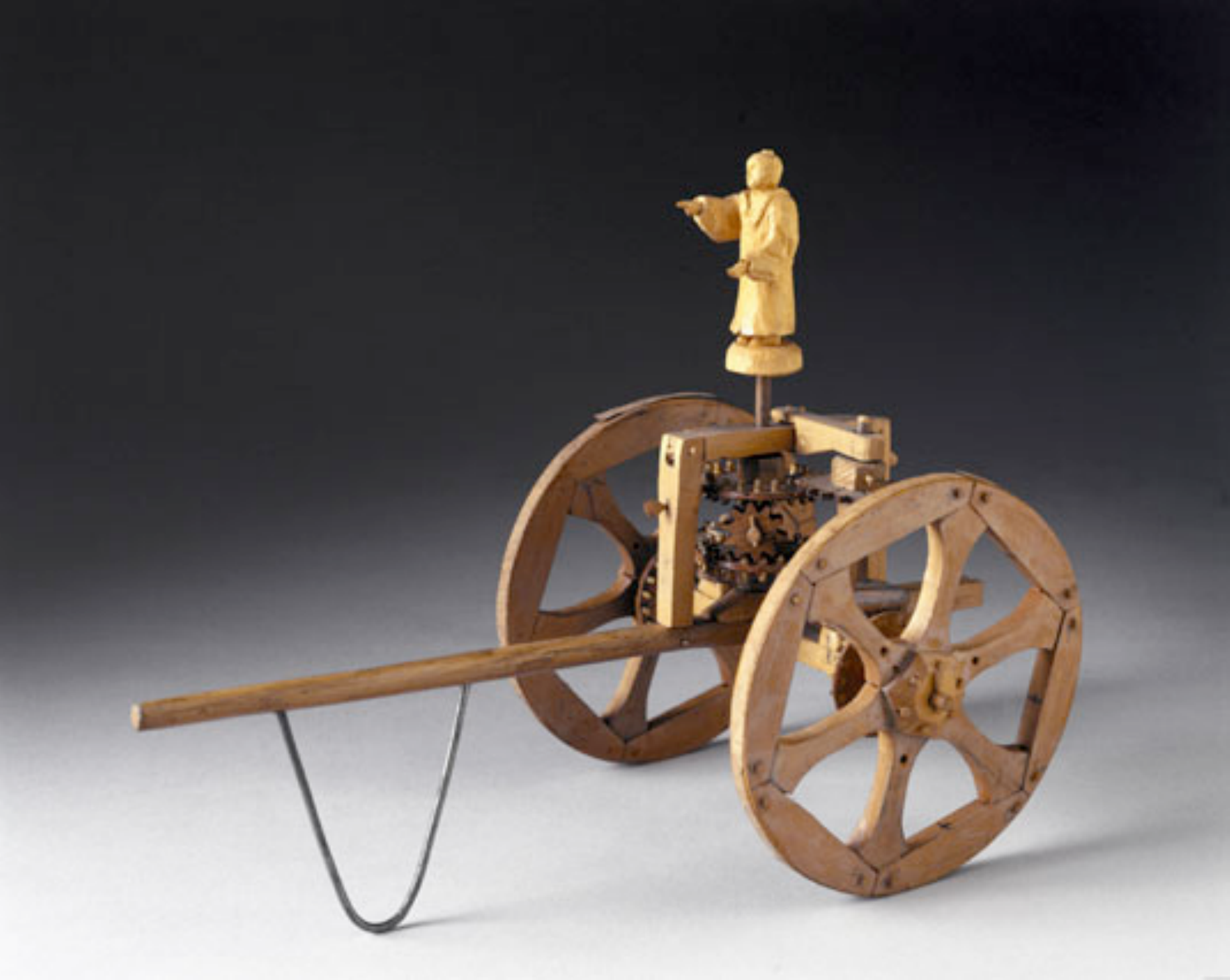} \hfill
\includegraphics[width=5cm]{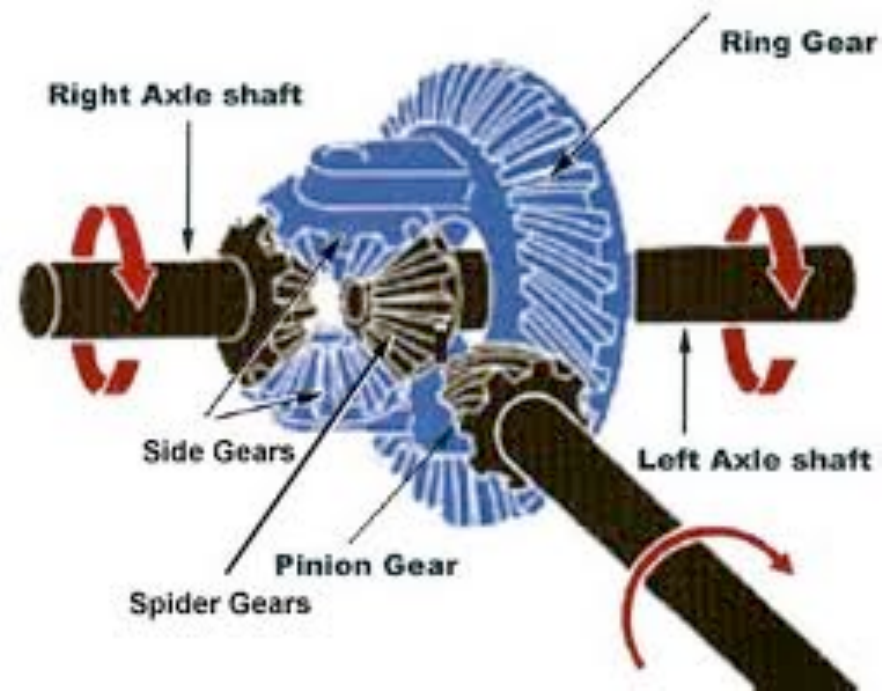}
\caption{South pointing chariot, and the differential that makes it work}
\label{fig:differential}
\end{figure}
Looking at the diagram in Figure~\ref{fig:differential}, 
we see that if the left and right axle are turning at the equal and opposite speed,
the ring gear and the central axle will not turn at all.  Believe me
that the rate of the central axle  depends linearly on
the rates of the left and right axles, and you must believe that
it turns at a rate proportional to the \emph{average} of the rates of
the left and right axles.   

While driving, \emph{when you turn right your left
  wheel travels further than your right wheel.}  Your car's differential
allows the vehicle to distribute the rotation of the drive train to
the two wheels as the requirements of the road demand. 
But Ma Jun was even cleverer with his supposed differential, because his chariot  not
only turned corners, but each time it did the figure on top
(traditionally a statue of an authoritative gentleman pointing south) turned
by exactly the same amount relative to the chariot in the opposite
direction, so that (from a fixed perspective) it remained pointing in
the same direction.  Look a little more carefully at
what happens when a vehicle turns a corner.

\begin{figure} 
\begin{tabular}{cc}
\begin{tikzpicture} 
\begin{scope}[>=latex]
\draw[style=thick] (-2,0) -- (-2,1) -- (-2.5,1);
\draw[style=thick] (2,0) -- (2,1) -- (3.2,1);

\path[fill=gray!20] (-2,0) -- (-2,1) -- (-2.5,1) -- (-2.5,3) -- (3.2,3) -- (3.3,1) -- (2,1) --
(2,0) -- (0,0);
\draw[white,style=ultra thick,dashed] (0,0) -- (0,3);
\draw[white,,style=ultra thick, dashed] (-2.5,3) -- (3.2,3);
\draw[color=blue,style= very thick] (.5,0) -- (.5,1);
\draw[color=blue,style= very thick]  (.5,1)  arc
  (180:90:1.5); 
\draw[color=blue,style= very thick]  (2,2.5)-- (3,2.5);
\draw[color=red,style= very thick] (1.5,0) -- (1.5,1);
\draw[color=red,style= very thick]  (1.5,1)  arc
  (180:90:.5); 
\draw[color=red,style= very thick]  (2,1.5)-- (3,1.5);

\draw[<->]   (.5,.5) -- (1.5,.5) node[pos=.5,above] {$w$};
\draw[color=blue,->] (2,1) -- +(130:1.5)
  node[above right,pos=.5,color=blue]{ {$r+w$}}; 
\draw[color=red,->] (2,1) -- +(170:.5)
  node[ below,pos=.5, color=red] { {$r$}};
\draw[color=blue,line width=2pt]  (.5,1)  arc
  (180:90:1.5); \node[color=blue] at (.2,2.3) {$(r+w)\theta$};  
\draw[color=red,line width=2pt]  (1.5,1)  arc
  (180:90:.5); 
\node[color=red] at (1.2,1.2) {$r\theta$};  
\end{scope}
 \end{tikzpicture}
&
\begin{tikzpicture} 
\begin{scope}[>=latex]
\draw[style=thick] (0,1.5) to [out=45, in = 210] (4,2.5);
\draw[color=blue,style=thick] (-.068,1.568) to [out=45, in = 210] (3.95,2.568);
\draw[color=red,style=thick] (.068,1.433) to [out=45, in = 210]
(4.05,2.42);
\draw[style=thick] (0,0) -- (.125,.125) arc (135:97:1.55) --
(2.75,.65) arc (-87.5:-60:2.8);
\draw[color=blue,style=thick] (-.068,.068) -- (.058,.191) arc (135:97:1.66) --
(2.75,.75) arc (-87.5:-60:2.7);
\draw[color=red,style=thick] (.068,-.068) -- (.192,.058) arc (135:97:1.49) --
(2.75,.55) arc (-87.5:-60:2.9);
\draw[color=blue,line width=2pt] (.058,.191) arc (135:97:1.66);
\draw[color=red,line width=2pt] (.192,.058) arc (135:97:1.49);
\draw[color=blue,line width=2pt] (2.75,.75) arc (-87.5:-60:2.7);
\draw[color=red,line width=2pt] (2.75,.55) arc (-87.5:-60:2.9);
\node at (2.05,1) {$+$};
\node at (.33, .85)  {$w\theta_1$};
\node at (3.15, 1.2)  {$w\theta_2$};
\node at (1,1.4) {${\color{blue} d_l} - {\color{red} d_r} = $};

\end{scope}
 \end{tikzpicture}
\end{tabular}
\caption{Relating the relative distance each wheeel travels to total rotation}
\label{fig:rotate}
\end{figure}
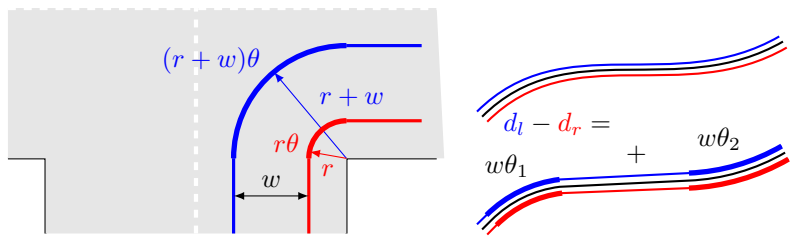

If as in the left side of Figure~\ref{fig:rotate} the chariot has two wheels a distance $w$
apart and it 
rotates right on the arc of a circle of radius $r$ for an angle $\theta$ (in radians!) the left wheel travels a distance $w
\theta$ more than the right wheel.  If the vehicle is not traveling in
an arc of a circle but along some arbitrary smooth path then trust
your intuition that the path can be approximated as closely as you like
by a sequence of straight lines and arcs of circles as in the right
side of Figure~\ref{fig:rotate}.  The left wheel travels $w\theta$
more than the right wheel on each arc segment, so adding these up with
signs,  \emph{whatever path your vehicle travels, if its
  total rotation is $\theta$ radians clockwise and its width is $w$
  the left wheel will travel a total of $w\theta$ more than the right
  wheel.}

Returning to our friend Ma Jun, attach the left axle of
the differential to the left wheel by an \emph{odd} number of gears,
 the right axle to the right wheel by an \emph{even} number of
gears, and finally connect the central axle to the gentlemanly
statue.  If the left wheel travels with velocity $v_l,$ its axle
 rotates at a rate proportional to $v_l.$  Each gear
reverses the direction of rotation, so the left axle of
the differential rotates at a rate $d\theta_l / dt$ proportional to
$-v_l.$  Likewise the right axle rotates at a rate $ d \theta_r/dt$
proportional to $v_r.$   That means that the statue rotates at a rate 
$d\theta_s/dt$ 
proportional to the average, which is to say proportional to $v_r -
v_l.$  Integrating over time the
statue rotates an angle $\theta_s$ proportional to $d_r-d_l,$ the
difference in the distances the two wheels travel over the course of
the journey.  adjust the
sizes of the various gears to make that constant of proportionality
$1/w,$ so that
the statute rotates an angle
\[\theta_s= \frac{d_r -d_l}{w}= - \theta,\]
the negative of the angle the chariot rotated.  That is to say, the
statue points in a constant direction.

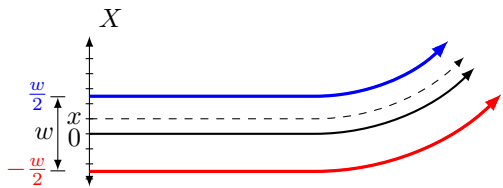
\begin{figure}
\begin{tikzpicture} 

\begin{scope}[>=latex]
\draw[|<->|] (-.42,0) -- (-.42,1);
\node at (-.6,.5) {$w$};
\draw[style=thick,->] (0,.5) -- (3,.5) arc (-90:-45:3);
\draw[color=blue,style=thick,->] (0,1) -- (3,1) arc (-90:-45:2.5);
\draw[color=red,style=thick,->] (0,0) -- (3,0) arc (-90:-45:3.5);
\draw[<->] (0,-.2) -- (0,1.8) node[above right] {$X$};
 \foreach \x in {-.1, .1,...,1.7} 
     \draw (-.05,\x) -- (.05,\x);
\node at (-.2,.45) {$0$};
\draw[style=dashed,->] (0,.7) -- (3,.7) arc
  (-90:-45:2.8);
 \node at (-.2,.7) {$x$};

\draw[color=blue,style=very thick,->] (0,1) -- (3,1) arc
(-90:-45:2.5);

\node[left,color=blue] at (-.42,1) {$\frac{w}{2}$};
\node[left,color=red] at (-.42,0) {$-\frac{w}{2}$};

\draw[color=red,style=very thick,->] (0,0) -- (3,0) arc (-90:-45:3.5);

\end{scope}
 \end{tikzpicture}
\caption{Rotation determines how moving the path slightly to the left
  changes the length}
  \label{fig:perturb}
\end{figure}

This is such a beautiful fact that as mathematics it deserves  to be 
expressed  in very  precise and abstract language so as to wring
every bit of meaning from it. 
Imagine an axis that cuts horizontally across the chariot, centered at
its center, so that the left and right wheel are placed at the points
$x=-w/2$ and $x=w/2$ respectively.  Choose some path that you would
like the chariot to travel, and let $d(x)$ be the distance a wheel placed at
position $x$ would travel during this journey.  So for instance
$d(-w/2)=d_l$ and $d(w/2)=d_r,$ as in Figure~\ref{fig:perturb}.  Now recognize the
formula for the net angle of rotation of our statue, $\sbrace{d(w/2)
  - d(-w/2)}/w,$ as a difference quotient.  The angle rotated, which is
the negative of the angle rotated by the chariot itself, is
independent of $w,$ so it is very natural to anyone who has taken
calculus to observe
\[\theta_s=\lim_{w \to 0} \frac{d(w/2)-d(-w/2)}{w}= \frac{\delta
  d}{\delta w},\]
That is \emph{the angle the statue rotates is the derivative of the
  distance traveled by a wheel with respect to its position on the
  chariot.}  I used $\delta$ for my derivative notation principally to
avoid the notational mostrosity $\frac{dd}{dx},$ but also as a sly
reference to functional differentiation and the calculus of variations, which are hiding in the wings
of this calculation and for which I hope to whet your appetite \cite{Weinstock74,Figueroa}. To put
it a different way (and to switch to the rotation of the vehicle) the
clockwise angle a vehicle rotates as it traverses a path is the rate
at which the length of the path increases as you move it to the left.

\section{Why It Doesn't Work}  \label{}

If the wheels and gears of the device are perfectly sized and suffer
from no slippage and the surface on which it travels is perfectly
flat, then a south pointing chariot which starts out
pointing south will always point south no matter what path it travels.  Slippage and errors in sizing are
not the concern of (theoretical) mathematicians, but nonflat surfaces
definitely are, and here is where the chariot transforms from a historical
and engineering curiosity to  mathematical tool.  

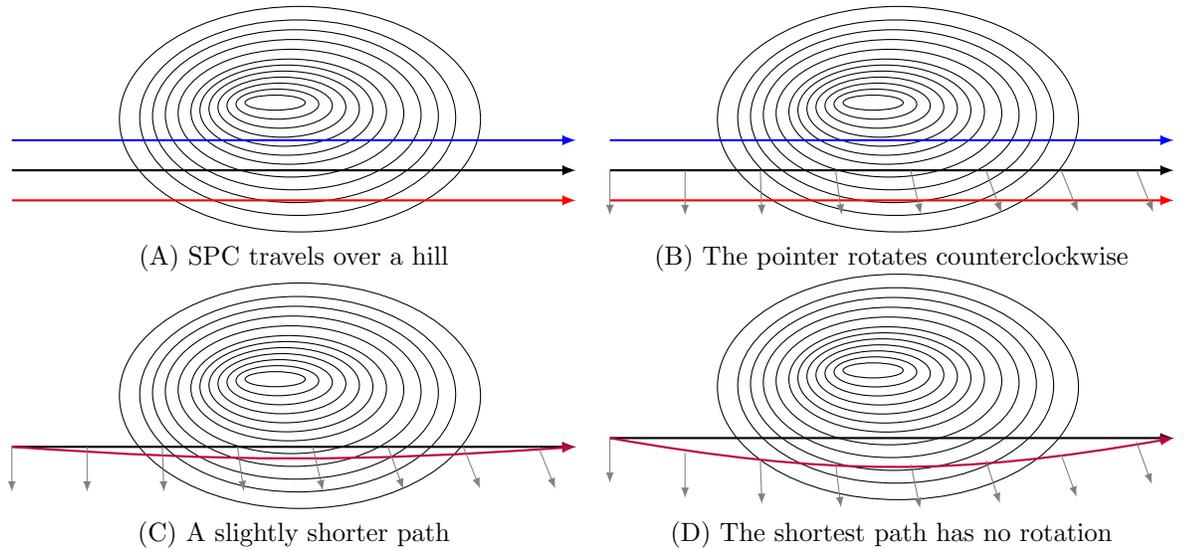
\begin{figure}
\begin{tabular}{cc}
\begin{tikzpicture} 
\begin{scope}[>=latex]

\draw (0,0) ellipse (.4 and .1);
\draw (.03,-.02) ellipse (.55 and .2);
\draw (.06,-.04) ellipse (.7 and .3);
\draw (.09,-.06) ellipse (.85 and .4);
\draw (.12,-.08) ellipse (1 and .5);
\draw (.15,-.1) ellipse (1.15 and .6);
\draw (.18,-.12) ellipse (1.3 and .7);
\draw (.21,-.14) ellipse (1.5 and .85);
\draw (.24,-.16) ellipse (1.7 and 1);
\draw (.27,-.18) ellipse (1.9 and 1.15);
\draw (.30,-.20) ellipse (2.1 and 1.3);
\draw (.33,-.22) ellipse (2.4 and 1.5);

\draw[style=thick,->] (-3.5,-.9) to (4,-.9);

\draw[style=thick,color=blue,->] (-3.5,-.5) to (4,-.5);
\draw[style=thick,color=red,->] (-3.5,-1.3) to (4,-1.3);

\end{scope}
 \end{tikzpicture}
&
\begin{tikzpicture} 
\begin{scope}[>=latex]
\draw (0,0) ellipse (.4 and .1);
\draw (.03,-.02) ellipse (.55 and .2);
\draw (.06,-.04) ellipse (.7 and .3);
\draw (.09,-.06) ellipse (.85 and .4);
\draw (.12,-.08) ellipse (1 and .5);
\draw (.15,-.1) ellipse (1.15 and .6);
\draw (.18,-.12) ellipse (1.3 and .7);
\draw (.21,-.14) ellipse (1.5 and .85);
\draw (.24,-.16) ellipse (1.7 and 1);
\draw (.27,-.18) ellipse (1.9 and 1.15);
\draw (.30,-.20) ellipse (2.1 and 1.3);
\draw (.33,-.22) ellipse (2.4 and 1.5);

\draw[style=thick,->] (-3.5,-.9) to (4,-.9);
\draw[style=thick,color=blue,->] (-3.5,-.5) to (4,-.5);
\draw[style=thick,color=red,->] (-3.5,-1.3) to (4,-1.3);
\draw[color=gray,->] (-3.5,-.9) -- +(-90:.6);
\draw[color=gray,->] (-2.5,-.9) -- +(-90:.6);
\draw[color=gray,->] (-1.5,-.9) -- +(-88:.6);
\draw[color=gray,->] (-.5,-.9) -- +(-82:.6);
\draw[color=gray,->] (.5,-.9) -- +(-77:.6);
\draw[color=gray,->] (1.5,-.9) -- +(-70:.6);
\draw[color=gray,->] (2.5,-.9) -- +(-68:.6);
\draw[color=gray,->] (3.5,-.9) -- +(-68:.6);
\end{scope}
 \end{tikzpicture}\\

(A) SPC travels over a hill  & (B) The pointer rotates
counterclockwise\\

\begin{tikzpicture} 
\begin{scope}[>=latex]
\draw (0,0) ellipse (.4 and .1);
\draw (.03,-.02) ellipse (.55 and .2);
\draw (.06,-.04) ellipse (.7 and .3);
\draw (.09,-.06) ellipse (.85 and .4);
\draw (.12,-.08) ellipse (1 and .5);
\draw (.15,-.1) ellipse (1.15 and .6);
\draw (.18,-.12) ellipse (1.3 and .7);
\draw (.21,-.14) ellipse (1.5 and .85);
\draw (.24,-.16) ellipse (1.7 and 1);
\draw (.27,-.18) ellipse (1.9 and 1.15);
\draw (.30,-.20) ellipse (2.1 and 1.3);
\draw (.33,-.22) ellipse (2.4 and 1.5);

\draw[style=thick,->] (-3.5,-.9) to (4,-.9);

\draw[style=thick,color=purple,->] (-3.5,-.9) to [out=-4, in=184 ] (4,-.9);
\draw[color=gray,->] (-3.5,-.9) -- +(-90:.6);
\draw[color=gray,->] (-2.5,-.9) -- +(-90:.6);
\draw[color=gray,->] (-1.5,-.9) -- +(-88:.6);
\draw[color=gray,->] (-.5,-.9) -- +(-82:.6);
\draw[color=gray,->] (.5,-.9) -- +(-77:.6);
\draw[color=gray,->] (1.5,-.9) -- +(-70:.6);
\draw[color=gray,->] (2.5,-.9) -- +(-68:.6);
\draw[color=gray,->] (3.5,-.9) -- +(-68:.6);
\end{scope}
 \end{tikzpicture}
& 
\begin{tikzpicture} 

\begin{scope}[>=latex]
\draw (0,0) ellipse (.4 and .1);
\draw (.03,-.02) ellipse (.55 and .2);
\draw (.06,-.04) ellipse (.7 and .3);
\draw (.09,-.06) ellipse (.85 and .4);
\draw (.12,-.08) ellipse (1 and .5);
\draw (.15,-.1) ellipse (1.15 and .6);
\draw (.18,-.12) ellipse (1.3 and .7);
\draw (.21,-.14) ellipse (1.5 and .85);
\draw (.24,-.16) ellipse (1.7 and 1);
\draw (.27,-.18) ellipse (1.9 and 1.15);
\draw (.30,-.20) ellipse (2.1 and 1.3);
\draw (.33,-.22) ellipse (2.4 and 1.5);

\draw[style=thick,->] (-3.5,-.9) to (4,-.9);

\draw[style=thick,color=purple,->] (-3.5,-.9) to [out=-10, in=190 ] (4,-.9);

\draw[color=gray,->] (-3.5,-.9) -- +(-90:.6);
\draw[color=gray,->] (-2.5,-1.1) -- +(-90:.6);
\draw[color=gray,->] (-1.5,-1.18) -- +(-88:.6);
\draw[color=gray,->] (-.5,-1.23) -- +(-83:.6);
\draw[color=gray,->] (.5,-1.25) -- +(-78:.6);
\draw[color=gray,->] (1.5,-1.2) -- +(-72:.6);
\draw[color=gray,->] (2.5,-1.12) -- +(-70:.6);
\draw[color=gray,->] (3.5,-.95) -- +(-70:.6);

\end{scope}
 \end{tikzpicture}\\
 (C) A slightly shorter path & (D) The shortest path has no rotation
\end{tabular}
\caption{Traversing a Hill}
\label{fig:hill}
\end{figure}

Imagine the chariot travels in a straight line (from a birds-eye view)
over a hill, keeping the peak of the hill to its left (see Figure~\ref{fig:hill}).  Notice that
the left wheel goes a little further uphill and a little further
downhill than the right wheel.  So the left wheel travels slightly
further than the right, and the statue rotates
slightly counterclockwise.  The bird flying overhead would assert that
the chariot traveled in a straight line, but a rider on the chariot,
watching the statue turn, would think that the chariot had turned
slightly to the right.  Perhaps you share the bird's point of view
(understandably) but suspend your preconceptions and
 judge the two points of view objectively.  In a disagreement
about what constitutes a straight line, what might we take as a
definition to adjudicate?  Euclid
famously said ``a straight line is the shortest distance between two
points,''  raising the question: Is the bird-straight
line the shortest distance between two points?

Were it the shortest distance between two points,  moving it a
bit to the left or right could not shorten it.  But SPC, the south pointing chariot,
measures the rate at which the length increases when moving the path to the
left.  Since the statue rotated counterclockwise, the path 
lengthens as it moves left and shortens as it moves right. If moving the bird-straight
line to the right shortens it, it is not the
shortest distance!

Well, hold on.  You might reasonably argue that we are looking for the
shortest distance between two fixed points and moving the path to the
right, shortens the distance, \emph{but also changes the endpoints.}
But that is really only a technical issue.
Consider a point $\epsilon$ to the right of the midpoint of SPCs
path.  Imagine SPC starting at the same starting point, but traveling
in a bird-straight path to this slightly moved point, then a
bird-straight path to its ending point. On the one hand the path has moved to the
right so it is  shorter, but on the other it is traversing two
sides of the triangle rather than the straight path between them, so  it
longer.  Which wins?  A little unfastidious thinking should convince
you that the shortening is \emph{linear} in $\epsilon$ (the factor is in fact
$\theta$) while the increasing is \emph{quadratic} in $\epsilon$ (by
Pythagoras).  So if $\epsilon$ is small enough the shortening wins: 
 Moving the path a little to the right
reduces is length.

If SPC travels
along this new path, it rotates a little less than it did on the
original path, but if the change is small enough it will still
rotate.  so moving it further to the right shortens it even more.
Keep doing that until the statue
does not rotate at all.  At that point the rotation $\delta d/\delta x=0,$ which is
to say that moving neither left nor right shortens the path.  If there is no way to move a path slightly while keeping the
endpoints fixed and make the length shorter, then the path is, at
least locally, the shortest distance between those two points.  In other
words a line which is straight according to the birds is not straight
in the sense of Euclid.   But straight according to the chariot
\emph{is} Euclid straight.  Euclid gives the nod to the chariot!

Geometers refer to a path which SPC can travel along with no rotation
of the statue by the lovely term ``geodesic.'' So a
path which is the shortest distance between two points is always a
geodesic (spoiler alert: not all geodesics are the shortest distance).   

A couple of details have been swept under the rug.
First, recall that $\delta d/
\delta x$ equals  the small $w$ limit of
the difference quotient $\sbrace{d(w/2) - d(-w/2)}/w,$ and on a
flat plane the difference quotient is constant and therefore equals
 the limit.  This constancy fails on a
curved surface --
for instance imagine SPC that is so wide that its wheels straddle the hill and miss it
entirely.  So from here on assume that SPC is ``sufficiently
small,'' that is so small compared to the topography it is
exploring that the difference quotient is extremely close to its
limit.  
Second, for the sake of simplicity I have spoken only of the total net
rotation that SPC undergoes along a path.  If, during 
some of the trip, the statue rotates clockwise while during other
parts it  rotates
counterclockwise, to find the minimum
would require perturbing to the left when it rotates clockwise and to
the right when it rotates
counterclockwise.  This process would converge on a path where
the statue does not rotate at all, not just one where the final net
rotation is zero.

One final aside. Following the direction of negative derivative to a
minimum where the derivative is zero should be familiar from
single and multivariable calculus.  This argument is 
just one more example of the familiar principle that the minimum of a
well-behaved function occurs where the derivative is zero, except with
one difference.  You are familiar with the principle when the function
depends on one variable, and if you know multivariable calculus
when it depends on two, or three, or many variables.  Here the function (the length of the path), depends on every point on
the path. Since there are infinitely many such points, one can say that
it depends on infinitely many variables.  To do this
honestly, you need to learn functional analysis.

\section{Exploring the Earth I}

\begin{figure} 
\begin{tabular}{cc}
\includegraphics[width=5cm]{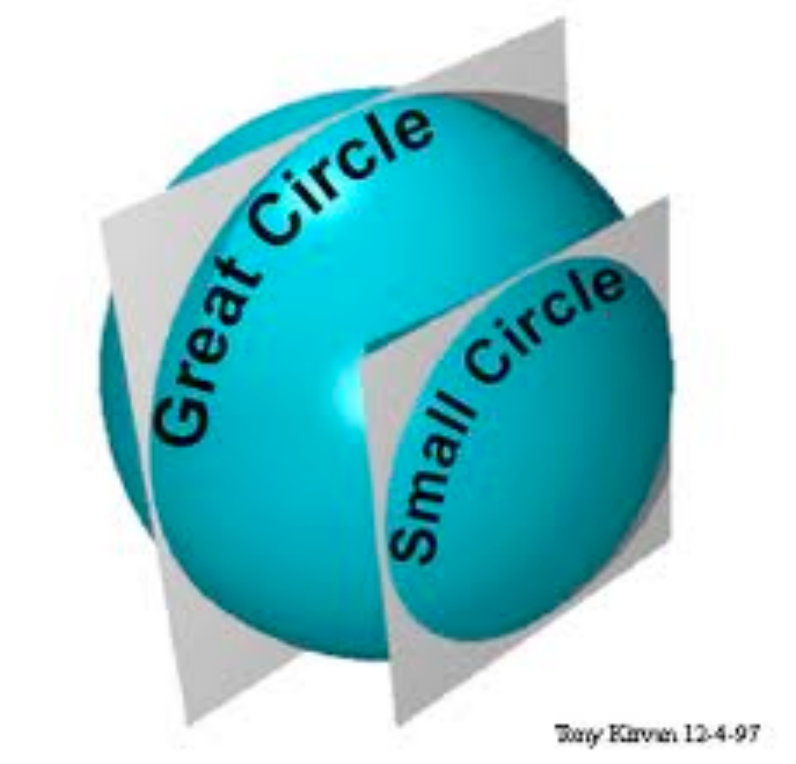}& \includegraphics[width=5cm]{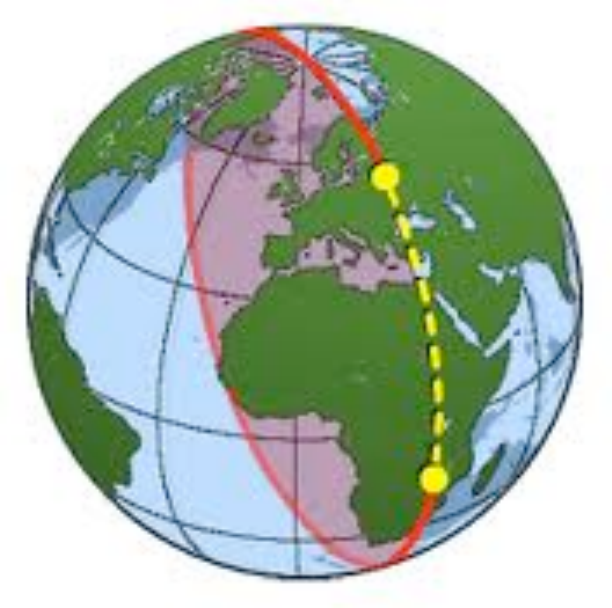}\\
 A great circle (geodesic) & A geodesic minimum and a geodesic saddle
 point
\end{tabular}
\caption{Geodesics on the sphere are arcs of great circles}
\label{fg:earth-geodesic}
\end{figure}
A sphere is one simple and interesting
example of a curved surface.  Suppose that SPC travels along the
equator of the earth.  This path is completely unchanged by
replacing the sphere with its mirror image around the equator,
so the left and right wheel travel the same distance and the statue
does not rotate.
The equator is a geodesic!  By the same reasoning, any circle
which cuts the sphere into two identical halves is a geodesic, as in
 Figure~\ref{fg:earth-geodesic}.  Such
circles, whose centers coincide with the center of the sphere, are
called ``great circles.''  Airplane pilots know this --
when they travel between two points far apart on the earth, they plot
a path that follows the great circle connecting them (again the
analogy with straight lines, and Euclid's remarkable prescience, are
affirmed:  Two points \emph{generally}  determine a great circle!), because they know
that this will be the shortest distance.   However, while traveling on a great
circle from Moscow to Dar es Salaam over the Mediterranean is a geodesic
and gives the shortest path between these two points,   traveling
the \emph{opposite way} on the same great circle (similar to the path
traveled in red in Figure~\ref{fg:earth-geodesic}) is still a geodesic but is
definitely \emph{not} the shortest route to Dar es Salaam from Moscow,
similar to the right side of Figure~\ref{fg:earth-geodesic}.
 Geodesics are not always minima of distance.  In fact, since
  such nonminimal geodesics can be perturbed in certain directions to
  make the distance smaller and in certain directions to make the
  distance larger, they are infinite-dimensional analogues of the
  saddle points of multivariable calculus.

Taking the analogy between lines on a plane and great
circles on a sphere in complete earnest, compare  Euclid's geometry of
triangles and parallel lines to their analogues on a sphere. (Some fine
print:  Two points on
a sphere determine a great circle 
\emph{unless they are antipodal.}   Students of spherical geometry
handle this embarrassment by treating pairs of antipodal points as a
single point, an issue we can safely ignore \cite{Mccleary13,Polking}) 
For example, suppose that SPC travels due west from some point on the
equator one quarter of the way around the earth, then due
north straight to the north pole, and then turns $\pi/2$ radians and
heads due south back to its starting point.  This
closed path made of three geodesics is called  a triangle in
spherical geometry.  But this triangle has three internal right
angles, something that could never happen in plane geometry, where the
sum of the internal angles of a triangle is always $\pi.$  In 
this case the sum is $3\pi/2,$ which is $\pi/2$ too many.  Euclid's
sum of angles theorem does not hold on the sphere -- it turns out that
this is the case because it relied
on the parallel line postulate.

If the statue starts
out pointing in the traditional south direction, it points
straight to the left on the first leg of the journey, directly behind
on the second, and 
straight to the right on the third leg, so it ends up pointing due
west.  This is more proof that SPC fails as a navigational aid -- not only does it fail to point
south under certain circumstances, it does not even return from a
journey pointing in the same direction as it started!  But the $\pi/2$
which it rotated in its journey and the $\pi/2$ by which the sum of
the angles exceeds the expected sum are related, perhaps not surprisingly.

Draw a few triangles on a sphere, make a rough estimate of
both the angle sum and the amount that SPC rotates as it traverses the
triangle,  and some guesses may present themselves.  The angle sum always
exceeds the expected $\pi,$ the amount it exceeds it by is always the
amount that the SPC rotates  (clockwise) over the journey, and this amount
seems to be roughly proportional to the size of the triangle.  If you
do not have a south pointing chariot on hand, you are just estimating the
quantity crudely, and I could not expect you to observe anything more
precise than ``the bigger the triangle the more the rotation.''  But
notice that a tall, very thin triangle has two almost
right angles and one almost $0$ angle, so this rotation is still quite
small.  That suggests that the right measure of the size of the triangle is
the area, not, say, its longest dimension.
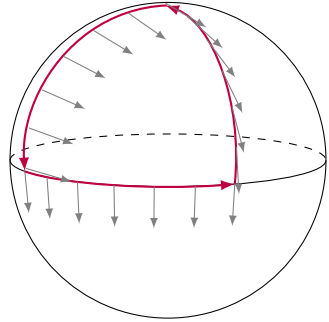
\begin{figure}
\begin{center}
\begin{tikzpicture}[scale=.7]

\begin{scope}[>=latex]
\draw (0,0) circle (3);
\draw (-3,0) arc (180:360: 3 and .5);
\draw[dashed] (-3,0) arc (180:0: 3 and .5);

\draw[thick,color=purple,->] (205:3 and .5) arc (205:295:3 and .5);

\foreach \x in {205,220,...,295}
\draw[color=gray,->] (\x:3 and .5) -- +(-.1*\x-64:.8);

\draw[thick,color=purple,->] (295:3 and .5) arc (-10:78:1.7 and 2.95);

\foreach \x in {3,18,...,78}
\draw[color=gray,->,xshift=-.4cm,yshift=.1cm] (\x:1.7 and 2.95) -- +(.75*\x-89:.9);

\draw[thick,color=purple,<-] (205:3 and .5) arc (184:92:2.8 and 2.95);

\foreach \x in {183,168,...,108}
\draw[color=gray,->,xshift=.1cm] (\x:2.8 and 2.95) -- +(.25*\x-62:.9);
\end{scope}
 \end{tikzpicture}
\end{center}
\caption{SPC gets very confused traversing this triangle}
\label{fig:triangle}
\end{figure}
This
observation is not unique to triangles.  The general statement on the
plane is that the sum of the \emph{external}
angles of a polygon  is always $2 \pi.$  On a sphere, the sum of the
external angles of a polygon (a closed curve made up of pieces of
great circles) will always be less than $2\pi,$ and less by the exact 
amount that the SPC rotates as it traverses that polygon, a quantity
that will in fact be exactly proportional to the area of the figure
encompassed \cite{Mccleary13,Polking}.  

All of this empirical exploration suggests that when SPC travels in a
loop, the total rotation it undergoes between the start and the finish
is some measure of how the surface curves in the interior of the
region subsumed by the loop.  This is correct, and the rest of this
article is devoted to convincing you that this is correct, and
exploring the outsized consequences of this one fact.

\section{Holonomy}

To understand what is going on requires a little abstraction.  For
each  loop $L$ (a path that begins and ends at the same point) on some surface,
let $H(L),$ the ``holonomy'' of the loop,  be the rotation that the pointer on SPC undergoes in traversing
that loop (make sure you are convinced that the amount that the pointer rotates is
not affected by where the pointer started). Holonomy views this
rotation as a function on
the set of all possible loops.

Here the narrative must pause again for a technicality.  The angle $2
\pi$ is the same angle as the angle $0.$ But if you have
ever walked a rambunctious dog on a leash, you know
that when the dog has run $2\pi$ radians around you, while she is in the
same position as if she had stayed still,  your leash can tell the
difference.  If you were riding on SPC, and you
slept through the journey, you would only know the angle of the final
position compared to the original angle.  But if you were awake and
watching the statue rotate (or tied a string to the statue's finger) you would know the
total rotation, and would be able to distinguish a $0$ rotation from a
$2\pi$ rotation from a $4\pi$ rotation.  Thus with care one can think
about the holonomy not just as an angle, but as a real number, with
differences of $2\pi$ being meaningful.  This distinction will only
matter near the very end, and can be ignored with little loss of comprehension.

To understand a function as a tool, a mathematician first asks how the
function responds to natural operations on its domain. One natural operation on loops is really an
operation  on the larger set of all paths.  If the path
$B$ begins at the point where the path $A$ ends, then it is natural to
define the \emph{composition} of the two paths $B \circ A,$ which traces
the path $A$ and then the path $B.$  If two loops $L_1$ and $L_2$ have the same
base point  (the common starting and ending point), then the composition
$L_1 \circ L_2$ is a loop  (you understand this definition if you can convince yourself that composition is
not commutative, i.e. that for some loops $L_1 \circ L_2 \neq L_2
\circ L_1$).  

So what is the holonomy of $L_2 \circ L_1$ in
terms of the holonomy of $L_1$ and the holonomy of $L_2?$  If you are
the kind of person who flips to the last page of a mystery novel
rather than trying to work it out yourself, resist the
desire and take a moment to think it through before proceeding.  Picture it rotating $30^\circ$ as it
traverses $L_1$ and  $20^\circ$ as it traverses $L_2$.  Then if
it traverses the first loop followed by the second it will rotate
... $50^\circ.$  In general 
\[H(L_2 \circ L_1) = H(L_2) + H(L_1).\]
Think of this property as saying that holonomy turns the composition operation on
loops into the addition operation on numbers.  When a function turns
the natural operation in the domain into the natural operation in the
range, it is called a \emph{homomorphism,} and this is one of the
nicest properties a function can have, according to a mathematical
aesthetic.

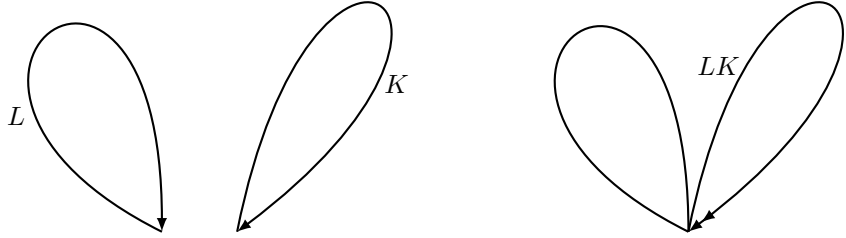
\begin{figure}
\begin{tikzpicture}
\begin{scope}[>=latex]
\draw[thick,->] (0,0) .. controls (-4,2) and (0,5) .. (0,0)
  node[near start, left] {$L$};
\draw[thick,->] (1,0) .. controls (2,5) and (5,3) ..  (1,0) node[near
end, right] {$K$};
\draw[thick] (7,0) .. controls (3,2) and (7,5) .. (7,0);
 \draw[thick,->>] (7,0) .. controls (8,5) and (11,3) ..  (7,0);
\node at (7.4,2.2) {$LK$};

\end{scope}
 \end{tikzpicture}
\caption{Composition of two loops} \label{fg:composition}
\end{figure}

The second property is a bit  subtler.  Suppose that SPC traverses
a loop $L,$ and somewhere in the middle of it turns off the loop,
travels along some path to a different destination, then returns on the
same path to where it left the loop, and continues the rest of the way
along the loop, as in Figure~\ref{fg:composition}.  Call this enhanced loop $L'$ and refer to it as
a \emph{detour.} How are $H(L)$ and $H(L')$ related?   When
SPC turns off its original path, it turns, say, $\theta$ radians clockwise, and
therefore the left wheel advances $w\theta/2$ and the right wheel $-w
\theta/2.$  Then SPC travels along the detour,and the left and right
wheel each advance by some unknowable amount.  At the terminus of the
detour, SPC turns around $\pi$ radians, lets say counterclockwise, so
the right wheel advances $w \pi /2$ and the left wheel advances $-w
\pi/2.$    Returning on the detour path in the opposite direction,
notice that the left wheel now travels exactly the distance the right
wheel had on the outbound journey, and vice versa.  So the total
distance traveled by each wheel from these two legs is exactly the
same, and thus has no effect on the rotation of the statue.  Finally,
upon returning to the original loop, SPC must rotate $\pi- \theta$
clockwise to return to its original direction.  Adding up all the left
and right wheel distances, you see that they are equal, and SPC has
not rotated a bit for all that effort.  From this conclude 
\begin{equation} \label{eq:detour}
H(L)= H(L'),
\end{equation}
that is to say adding a detour in this sense to any loop has no effect
on the holonomy.

\begin{figure}
\begin{tabular}{ccc}
\begin{tikzpicture}[scale=0.6]
\begin{scope}[>=latex]
\draw[thick,->] (0,0) to [out=90, in=225] (1,3);
\draw[thick] (0,0) to [out=90, in=225] (1,3);
\draw[thick,->] (1,3) to [out=135, in=45] (-4,3);
\draw[thick] (0,0) to [out=90, in=225] (1,3);
\draw[very thick,->] (-4,3) to [out=45, in=135] (1,3);
\draw[thick] (0,0) to [out=90, in=225] (1,3);
\draw[very thick] (-4,3) to [out=45, in=135] (1,3);
\draw[thick,->] (1,3) to [out=45, in=0] (0,0) node[near end, left] {$L$};
\end{scope}
 \end{tikzpicture}
&
\begin{tikzpicture}[scale=0.6]
\begin{scope}[>=latex]
\draw[thick] (0,0) to [out=90, in=225] (1,3);
\draw[very thick] (-4,3) to [out=45, in=135] (1,3);
\draw[thick,->] (1,3) to [out=45, in=0] (0,0) node[near end, left] {$L$};

\draw[color=blue,->] (-.3,0) to [out=90, in=225] (1-.212,3+.212);
\draw[color=red,->] (.3,0) to [out=90, in=225] (1+.212,3-.212);
\draw[color=blue] (-.3,0) to [out=90, in=225] (1-.212,3+.212);
\draw[color=red] (.3,0) to [out=90, in=225] (1+.212,3-.212);
\draw[color=blue,->] (1-.212,3+.212) arc (135:225:.3);
\draw[color=red,->] (1+.212,3-.212) arc (-45:45:.3);
\draw[color=blue] (-.3,0) to [out=90, in=225] (1-.212,3+.212);
\draw[color=red] (.3,0) to [out=90, in=225] (1+.212,3-.212);
\draw[color=blue] (1-.212,3+.212) arc (135:225:.3);
\draw[color=red] (1+.212,3-.212) arc (-45:45:.3);
\draw[color=blue,->] (1-.212,3-.212) to [out=135, in=45] (-4+.212,3-.212);
\draw[color=red,->] (1+.212,3+.212) to [out=135, in=45] (-4-.212,3+.212);
\end{scope}
 \end{tikzpicture}
&
\begin{tikzpicture}[scale=0.6]
\begin{scope}[>=latex]
\draw[thick] (0,0) to [out=90, in=225] (1,3);
\draw[very thick] (-4,3) to [out=45, in=135] (1,3);
\draw[thick,->] (1,3) to [out=45, in=0] (0,0) node[near end, left] {$L$};

\draw[color=blue,->] (-.3,0) to [out=90, in=225] (1-.212,3+.212);
\draw[color=red,->] (.3,0) to [out=90, in=225] (1+.212,3-.212);
\draw[color=blue] (-.3,0) to [out=90, in=225] (1-.212,3+.212);
\draw[color=red] (.3,0) to [out=90, in=225] (1+.212,3-.212);
\draw[color=blue,->] (1-.212,3+.212) arc (135:225:.3);
\draw[color=red,->] (1+.212,3-.212) arc (-45:45:.3);
\draw[color=blue] (-.3,0) to [out=90, in=225] (1-.212,3+.212);
\draw[color=red] (.3,0) to [out=90, in=225] (1+.212,3-.212);
\draw[color=blue] (1-.212,3+.212) arc (135:225:.3);
\draw[color=red] (1+.212,3-.212) arc (-45:45:.3);
\draw[color=blue,->] (1-.212,3-.212) to [out=135, in=45] (-4+.212,3-.212);
\draw[color=red,->] (1+.212,3+.212) to [out=135, in=45] (-4-.212,3+.212);
\draw[color=blue] (-.3,0) to [out=90, in=225] (1-.212,3+.212);
\draw[color=red] (.3,0) to [out=90, in=225] (1+.212,3-.212);
\draw[color=blue] (1-.212,3+.212) arc (135:225:.3);
\draw[color=red] (1+.212,3-.212) arc (-45:45:.3);
\draw[color=blue] (1-.212,3-.212) to [out=135, in=45] (-4+.212,3-.212);
\draw[color=red] (1+.212,3+.212) to [out=135, in=45] (-4-.212,3+.212);
\draw[color=blue,->] (-4+.212,3-.212) arc (315:135:.3);
\draw[color=red,->] (-4-.212,3+.212) arc (135:-45:.3);
\draw[color=blue] (-.3,0) to [out=90, in=225] (1-.212,3+.212);
\draw[color=red] (.3,0) to [out=90, in=225] (1+.212,3-.212);
\draw[color=blue] (1-.212,3+.212) arc (135:225:.3);
\draw[color=red] (1+.212,3-.212) arc (-45:45:.3);
\draw[color=blue] (-4+.212,3-.212) arc (315:135:.3);
\draw[color=red] (-4-.212,3+.212) arc (135:-45:.3);
\draw[color=blue!50!red,->] (-4+.212,3-.212) to [out=45, in=135] (1-.212,3-.212);
\draw[color=blue!50!red,->] (-4-.212,3+.212) to [out=45, in=135] (1+.212,3+.212);
\draw[color=blue] (-.3,0) to [out=90, in=225] (1-.212,3+.212);
\draw[color=red] (.3,0) to [out=90, in=225] (1+.212,3-.212);
\draw[color=blue] (1-.212,3+.212) arc (135:225:.3);
\draw[color=red] (1+.212,3-.212) arc (-45:45:.3);
\draw[color=blue] (-4+.212,3-.212) arc (315:135:.3);
\draw[color=red] (-4-.212,3+.212) arc (135:-45:.3);
\draw[color=blue!50!red] (-4+.212,3-.212) to [out=45, in=135] (1-.212,3-.212);
\draw[color=blue!50!red] (-4-.212,3+.212) to [out=45, in=135] (1+.212,3+.212);
\draw[color=blue,->] (1+.212,3+.212) arc (45:135:.3);
\draw[color=red,->] (1-.212,3-.212) arc (225:315:.3);
\draw[color=blue] (-.3,0) to [out=90, in=225] (1-.212,3+.212);
\draw[color=red] (.3,0) to [out=90, in=225] (1+.212,3-.212);
\draw[color=blue,->] (1-.212,3+.212) arc (135:225:.3);
\draw[color=red,->] (1+.212,3-.212) arc (-45:45:.3);
\draw[color=blue,->] (-4+.212,3-.212) arc (315:135:.3);
\draw[color=red,->] (-4-.212,3+.212) arc (135:-45:.3);
\draw[color=blue!50!red] (-4+.212,3-.212) to [out=45, in=135] (1-.212,3-.212);
\draw[color=blue!50!red] (-4-.212,3+.212) to [out=45, in=135] (1+.212,3+.212);
\draw[color=blue,->] (1+.212,3+.212) arc (45:135:.3);
\draw[color=red,->] (1-.212,3-.212) arc (225:315:.3);
\end{scope}
 \end{tikzpicture}
\end{tabular}
\caption{A loop with a detour, and the tracks of SPC} \label{fg:detour}
\end{figure}
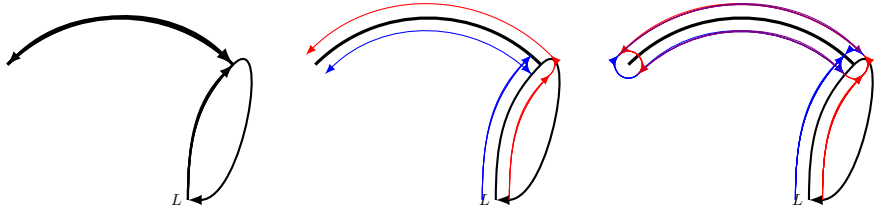

Perhaps you noticed a  subtlety in the above argument.  At the end
of the detour SPC turned around in the counterclockwise
direction, but might just as well have chosen to turn in the
clockwise direction.  Had it done so, a careful accounting shows
that the statue would have added one extra counterclockwse rotation of
$2 \pi$ during the detour. To compute only the angle of
holonomy, this is of no account, but in terms of the numerical
holonomy,  Eq.~\eqref{eq:detour} is  true only if 
SPC does its about-faces correctly.  Specifically,
whatever direction SPC turns to get on and off the detour, it must
turn the opposite way at the extreme point of its detour.  

There is a second case involving a detour move.  The stereotypical American way
to exercise is to drive to
the running track, run around the track, and then drive home.  Here
the loop $L$ is the run around the track, and the longer loop $L'$ is
the loop $L$ bracketed by a path and its reverse at the beginning and
end of the loop.  So now $L$ and $L'$ have different starting points,
but  repeating the logic it is still true that
$H(L')=H(L).$  Once again for numerical holonomy attention must be paid
to the direction of rotations.

So holonomy  assigns a number or angle to each loop in such a way that composition of
loops turns into addition and detours have no effect. These two
properties are already enough to make a key observation.  Imagine that
the loop SPC travels on a loop surrounding a region.  You may
not have imagined that there was any other possibility, but  if the surface were
the surface of a doughnut, for instance, what mathematicians call a
\emph{torus,} the loop could wrap around the doughnut as if it were a
ribbon wrapped around a gift.  Instead  assume that the loop  surrounds a region $R.$
First notice that the holonomy is the same  regardless of what point is the 
start and end of the loop, because the effect of starting and ending at a
different point is exactly a detour move of the second kind. So define the holonomy
of the \emph{region} $R$ to be the holonomy of any loop $L$ that goes once
clockwise around the region.  This does not specify the loop
uniquely, but does specify the holonomy.

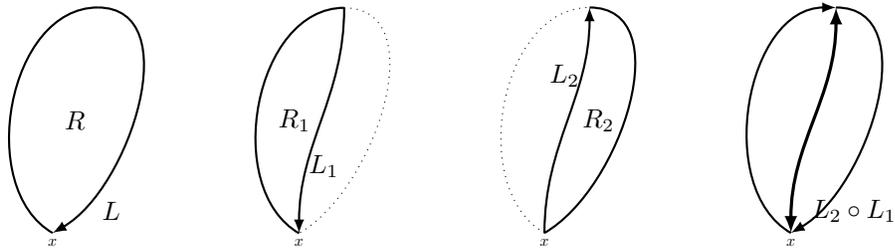
\begin{figure}
\begin{tabular}{cccc}
\begin{tikzpicture}[scale=.6]
\begin{scope}[>=latex]
\draw[thick] (0,0) to [out=150, in=180] (1,5);
\draw (1.3,.5) node {$L$};
\draw[thick,->] (1,5) [out=0, in=30] to (0,0)  node[near end, below]
{$x$};
\draw (.5,2.5) node {$R$};
\end{scope}
 \end{tikzpicture}
&
\begin{tikzpicture}[scale=.6]
\begin{scope}[>=latex]
\draw[thick] (0,0) to [out=150, in=180] (1,5);
\draw (.55,1.5) node {$L_1$};
\draw[dotted] (1,5) [out=0, in=30] to (0,0)  node[near end, below]
{$x$};
\draw (-.1,2.5) node {$R_1$};
\draw[thick,->] (1,5) [out=270, in=90] to (0,0);
\end{scope}
 \end{tikzpicture}
&
\begin{tikzpicture}[scale=.6]
\begin{scope}[>=latex]
\draw[dotted] (0,0) to [out=150, in=180] (1,5);
\draw (.45,3.5) node {$L_2$};
\draw[thick] (1,5) [out=0, in=30] to (0,0)  node[near end, below]
{$x$};
\draw (1.2,2.5) node {$R_2$};
\draw[thick,->] (0,0) [out=90, in=270] to (1,5);
\end{scope}
 \end{tikzpicture}
&
\begin{tikzpicture}[scale=.6]
\begin{scope}[>=latex]
\draw[thick,->] (0,0) to [out=150, in=180] (1,5);
\draw (1.4,.5) node {$L_2 \circ L_1$};
\draw[thick,->] (1,5) [out=0, in=30] to (0,0)  node[near end, below]
{$x$};
\draw[very thick,<->] (0,0) [out=90, in=270] to (1,5);
\end{scope}
 \end{tikzpicture}
\end{tabular}
\caption{Cutting a region into pieces}\label{fg:cutting}
\end{figure}

Now  divide $R$
into two regions $R_1$ and $R_2.$  Let $x$ be a point on their shared
boundary, let $L_1$ be the loop that goes clockwise
around  $R_1$ based at $x,$ and let $L_2$ be the loop that goes
clockwise around $R_2$ based at $x$ as in Figure~\ref{fg:cutting}.  Notice that $L_2\circ L_1$ is a
detour of $L.$  Thus $H(L)=H(L_1)+H(L_2)$ so  the holonomy of $R$ is the sum of the
holonomies of $R_1$ and $R_2.$  This argument generalizes
to say that the holonomy of any region divided into
subregions is the sum of the holonomies of the subregions.  Holonomy behaves like area,
in that the whole is equal to the sum of its parts.

The last thing I need you to believe about holonomy puts the most
demand on your visual intution.  The rough idea is that holonomy
should be continuous, that small changes to the loop should cause
small changes in the holonomy.  In particular if you move a small
loop from a flat area to a steep hill, you would not be surprised if
the holonomy changed a lot, but if you moved a small loop at one point
on a hill to a similarly shaped small loop at a nearby point, near
enough that the ``roundness'' or ``curvature'' of the hill has not
changed much, you would expect the holonomy to change very little.

Now consider four nearby points on the surface, arranged so as to form
an approximate rectangle, i.e. so that when they are connected by
geodesics all four angles are close to right angles.  
Let $R$ be the interior of this rectangle and consider the quantity
\[\frac{H(R)}{\text{Area}(R)}.\] 
By dividing two opposite edges into $n$ equal length intervals and the
other pair of opposite edges into $m$ equal intervals, we divide the
entire approximate rectangle into $nm$ approximate rectangles that are
all of approximately the same shape and size. So, it should be plausible that their holonomies are all
approximately equal and that their areas are all approximately equal.
Since both holonomies and areas add, calling one of the smaller
rectangles $R_1$ we have that 
\[\frac{H(R)}{\text{Area}(R)}\sim \frac{nmH(R_1)}{nm\text{Area}(R_1)}=
\frac{H(R_1)}{\text{Area}(R_1)}.\]
The approximation gets better the smaller the size of the rectangle, so
imagine a sequences of rectangles $R_n$, whose sizes approach zero and
which all surround one point $x.$ This heuristic argument suggests
that the limit of the ratio of holonomy to area converges for this
sequence, which in fact it does.  Define the \emph{curvature} $K(x)$
of the surface at the point $x$ to be this limit
\[\lim_{n \to \infty} \frac{H(R_n)}{\text{Area}(R_n)}= K(x).\]
In fact, by varying $n$ and $m$ it should also be plausible that this
number does not depend on the shape of the sequence of rectangles.  

\begin{figure}
\begin{tikzpicture}
\begin{scope}[>=latex]

\draw[very thick,->] (0,0) .. controls (4,4) and (-3,5) .. (0,0) node[near end, left] {$L$}; 
\begin{scope}
\clip (0,0) .. controls (4,4) and (-3,5) .. (0,0);
\draw[step=2mm] (-3,0) grid (4,5);
\end{scope}
\draw[thick] (.2,2) rectangle (.4,2.2);
\draw[->] (.2,2) -- (6,1);
\draw[->] (.4,2.2) -- (8,3);
\draw[thick,color=purple] (6,1) rectangle (8,3) node[pos=0,left, below] {$(x_i,y_i)$};
\draw[thick,->,color=purple] (6,3) -- (6,1);
\draw[|<->|] (6,3.3) -- (8,3.3) node[pos=.5,above] {$\Delta x$};
\draw[|<->|] (8.3,1) -- (8.3,3) node[pos=.5,right] {$\Delta y$};
\draw[color=purple,->] (4.8,3) .. controls (6,2) .. (7,2) node[pos=0,above,color=purple] {$K(x_i, y_i)$};
\end{scope}
 \end{tikzpicture}
\caption{Chopping into small rectangles}\label{fg:rectangle}
\end{figure}
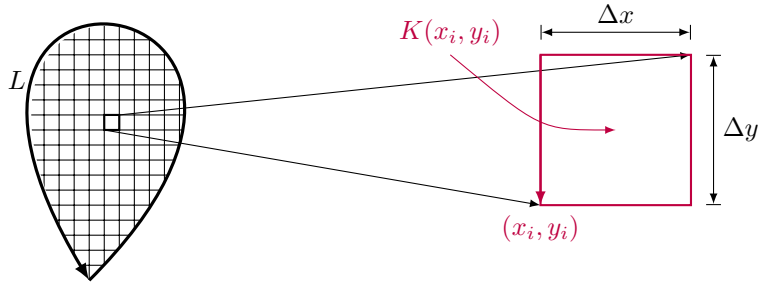

Below I will argue that this real-valued function $K(x)$  defined on the
surface  does a good job of
measuring your intuitive notion of how curved the surface is at each point.  But first a remarkable
fact:  Another way to express this limit is to say that for such
a rectangle the holonomy is approximately equal to the
area times the curvature of some point in the interior.  Imagine an
arbitrary region, but still small, so that we can draw a sequence of
geodesics on it that divide the region up into many small rectangles
(this requires replacing the boundary with a more jagged
approximation). For the $i$th approximate rectangle, call the side
lengths $\Delta x_i$ and $\Delta y_i,$ and let $K(x_i, y_i)$ be the
curvature of some point in the rectangle, so that the holonomy of that
rectangle should be approximately 
\[K(x_i,y_i) \Delta x_i \Delta y_i\]
as in Figure~\ref{fg:rectangle}
In particular the holonomy of the entire region should be approximated
by 
\[\sum_i K(x_i, y_i) \Delta x_i \Delta y_i.\]
Of course here when we say approximate we are including errors that
come from approximating the region by these (approximate) rectangles,
from approximating the area by $\Delta x_i \Delta y_i,$ and from
approximating the holonomy by area times the curvature at one point.
Your trust is required to believe that as the region is divided more
finely all these errors go to zero rapidly enough that 
\[H(R) = \lim_{\Delta x_i,\Delta y_i \to 0} \sum_i K(x_i, y_i) \Delta
  x_i \Delta y_i.\]

But  granting this much trust, then all that is required of
you is recollection.  You have seen (at least if you have taken
multivariable calculus before) this  limit.  It is exactly how
the integral of a function over a two dimensional region is defined.
The conclusion is that the holonomy around a closed loop that
surrounds a region $R$ is the integral of curvature over $R,$ or
\[H(R)=  \iint_{R} K(x,y) dx dy.\]
On the one hand, this is a really remarkable claim.  If 
SPC goes on a long trip around northern Arizona,
this says that the position of the flag at the end depends on the
value of  curvature at each point that its path surrounds.  If
the path encloses the Grand Canyon, where one presumes the curvature
function would vary wildly, the position of the flag will in some
sense record information about this curvature at each point in the
canyon, even though from the SPC one might never have
seen the Grand Canyon.  On the other hand, this result
looks superficially very similar to, and indeed is deeply connected
to, results like the Divergence Theorem and Green's Theorem, which
relate an integral along a closed path to a higher dimensional
integral of something else over the interior (note in paragraph 4 of
Section 2 that the rotation of the pointer is calculated by an integral).

Let's explore what this notion of curvature looks like with some
examples.

\section{Exploring the Earth II}

The curvature at a give point on the sphere is given by the limit  of the ratio of the holonomy
around the rectangle to the area of the rectangle as the size of the
rectangle goes to zero.  Consider such a decreasing sequence of smaller rectangles at
some point $p$ on the sphere, and suppose that $q$ is a different point on
the sphere.  Then there is a rotation of the sphere that takes $p$ to
$q,$ which therefore takes that sequence to a sequence of decreasing
rectangles around $q.$  Because a rotation of the sphere will clearly
preserve all distances, angles and areas, the resulting loops will still be
approximately rectangular and will have the same area and holonomy.
Thus the ratio will converge to the same thing, so that the curvature
at every point of the sphere must be the same.  This makes perfect sense:
The sphere is maximally symmetric, so a reasonable measure of the
curvature at each point should be preserved by all those symmetries.  

 In Section 4 the holonomy of the loop pictured in
Figure~\ref{fig:triangle} was computed to be $\pi/2.$  This holonomy
is equal to the integral of the curvature function over the interior
of the triangle $\int \int_R K(p) dp.$  However, since the curvature is a
constant $K,$ the holonomy is $K\text{Area}(R).$  The area of this region is  $1/8$ of the
area $4 \pi r^2$ of the sphere, from which it follows
that the curvature at each point of a sphere of radius $r$ is
$1/r^2.$  This very simple formula should make intuitive sense.  A
large radius sphere, like the earth, is almost flat, while a very
small radius sphere is extremely curved.

Similarly, a small rectangle on the cylinder with two sides parallel
to the axis of rotation and two sides perpendicular to it can be translated up and down the
axis and rotated around it to a rectangle with the same angles
anywhere on the cylinder.  So by the same reasoning, the cylinder has
constant curvature.  But in this case the small rectangle made up of
geodesics has four right angles, which means that its angle deficit is
zero and therefore the curvature of the cylinder is $0.$  The
cylinder, if you will, is flat.  This quite likely does not fit so
well with your intuitive notion of curvature.

In fact there is more than one intuitive notion of
curvature that we tend to mix up.  If you deform a surface so that all
the distances on it remain unchanged, you have not changed its
curvature in this sense.  To make this somewhat technical notion a
little more concrete, consider a piece of cloth.  Draw a
path on the piece of cloth and then bunch it up.  You've deformed the piece of cloth,
but the path length has not changed, because that would require
stretching the fabric.  If lengths of paths and therefore distances
are unchanged by anything you can do to (idealized) cloth, then the
amount SPC rotates along a path is also unchanged,
since it only depends on the distances each wheel travels.  Thus
holonomy and curvature are unchanged by such a deformation.  We speak
of ``intrinsic geometry'' to refer to any property of a surface that
is unchanged by a deformation you can do to cloth without stretching
it.  That is to say, a deformation that does not change the lengths of
any paths.  

A flat piece of cloth can be wrapped around a cylinder without
stretching any part of it.  So  a cylinder has the same
intrinsic geometry as a flat plane, and in particular has the same
curvature.  The kind of bending that turns the flat cloth into a
cylinder does not change its intrinsic geometry and ``does not count''
as far as curvature is concerned.

Which brings us at last to the point.

\section{Mapping the Earth}

\begin{figure}
\includegraphics[width=5cm]{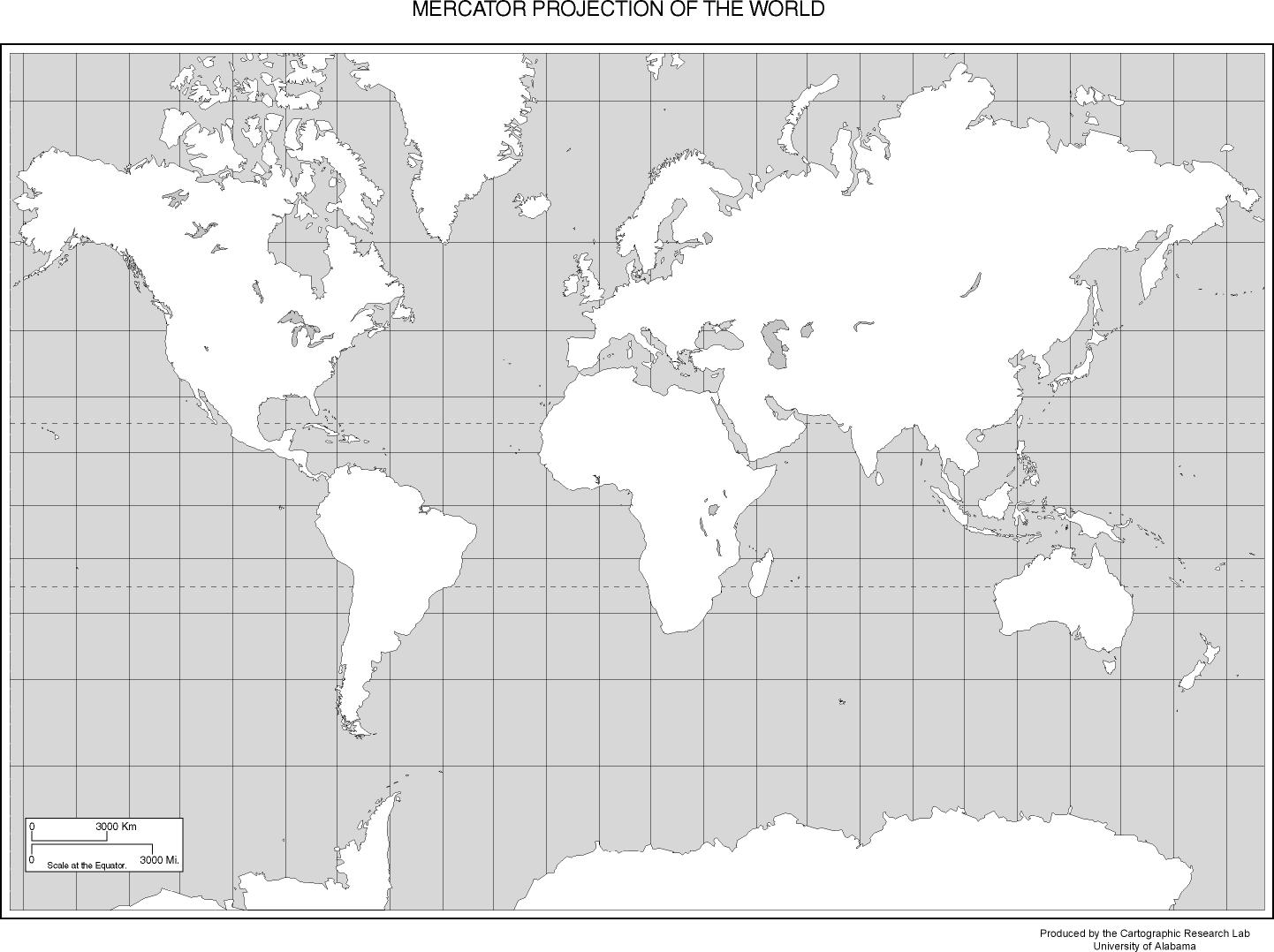} \hfill
\includegraphics[width=5cm]{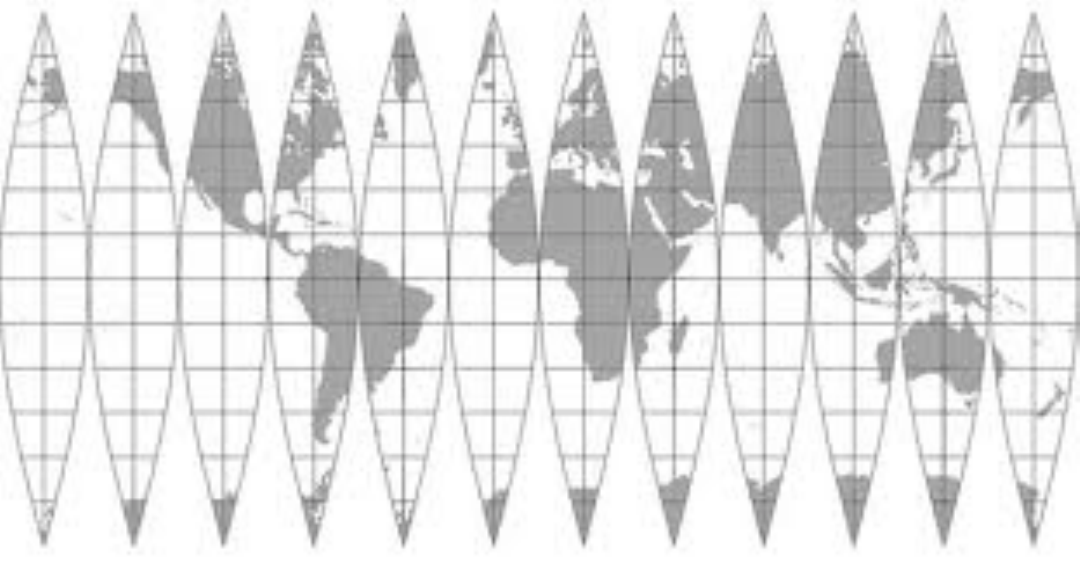}
\caption{Mercator and a competing projection of the earth}
\label{fig:projection}
\end{figure}

The standard Mercator projection (Figure~\ref{fig:projection}) does
not represent distances accurately.  It makes Canada, Russia, and
Antarctica look immense, while equally big regions near the equator look
small.  Obviously this is a poor property for a map to have.  Various
schemes attempt to solve this problem, but can they?  Is it possible
to draw a map of the world on which all distances are represented
proportionally?  Another form of the question asks if one can peel an
orange so that the skin can be laid flat without stretching. Note that
you
are allowed to cut the skin in various places, as in
Figure~\ref{fig:projection}. Nevertheless the  answer is no.  
The reason is clear from our previous discussion.  If the
identification between the earth and the map preserved distances, then
since the workings of SPC only depend only upon distances, every loop on
the earth must have the same holonomy as the corresponding loop on the
map.  But the holonomy of a loop on the earth will be proportional to
the area enclosed, whereas the holonomy of a loop on a flat map is
necessarily zero.  So such a map is impossible! 

\section{Gauss' Theorema Egregium and Beyond}

This notion of curvature, and the many things it tells you, including
the fact about maps of the earth, are due to Gauss.  Gauss did not
describe curvature in terms of south pointing chariot.  Gauss gave a
definition of the curvature that depends on how the surface is
embedded in three space.  Roughly Gauss imagines rotating and moving
the surface in three dimensions until the desired point on the surface
is at the origin and the tangent plane at that point is the $xy$
plane.  At least near the origin the surface can then be written in
the form $z=f(x,y),$ where the Taylor series expansion of $f$ has no
constant or linear terms.  A further rotation about the $z$ axis
guarantees that  
\[f(x,y)=ax^2 + by^2 + \, \text{third order and higher}.\]
Gauss' curvature is simply $4ab.$  Gauss was then able to prove that
although this definition clearly depends on the details of how the
surface is embedded in three space, the resulting curvature does not,
and only depends on the distances within the surface, that is on its
intrinsic geometry \cite{Mccleary13}.  Gauss found this result so surprising that he
called it his ``Theorema Egregium'' or ``Remarkable Theorem.''  If you
are familiar with even a few of the revolutionary results which are
due to Gauss (several of them simply called Gauss' Lemma), you will
understand that this is quite a statement.  

Our definition of curvature is based on the SPC which is manifestly
intrinsic, and thus seems to have  avoided the work of the Theorema
Egregium.  Of course this is because we skipped over hard analytic
issues, and in particular relied on nothing but the reader's trust for
the existence of the limit (and rate of convergence) of the definition
of curvature.   

Gauss' work has been extended in several directions.  Most
importantly, all this can be done with great care in more dimensions
than two (replace SPC with a rocket ship with a gyroscope), where the
notion of curvature and holonomy are quite a bit more complicated (one
must keep track of the change in the path's length as you perturb it
in each direction).  SPCs doomed attempt to keep track of a global
direction (towards the south pole) with a local calculation (the gears
and the distances its wheels travel) is using the mathematical notion
of parallel transport \cite{Santander92}, which can be generalized beyond distances to
mathematical objects called connections.   Connections are used to
model all of
the fundamental forces of nature except gravity (it plays a role in
modeling gravity as well) and, in some reckonings, even monetary inflation!
These ideas are the building blocks of differential geometry.

\bibliographystyle{alpha}\def\cprime{$'$}

\end{document}